\documentclass{amsart} 
\input diagrams
\numberwithin{equation}{section}
\newtheorem{theorem}[equation]{Theorem}
\newtheorem{proposition}[equation]{Proposition}

\theoremstyle{remark}
\newtheorem{remark}[equation]{Remark}
\newtheorem{definition}[equation]{Definition}
\newtheorem*{remark*}{Remark}
\newtheorem*{note}{Note}

\newarrow{Mapsto}{|}{-}{-}{-}{>}
\newarrow{Equal}{=}{=}{=}{=}{=}

\begin{document}
\title[Twisted Fourier-Mukai transforms and non-K\"ahler elliptic surfaces]
{Twisted Fourier-Mukai transforms and bundles on non-K\"ahler elliptic surfaces}

\author{Vasile Br\^{\i}nz\u{a}nescu}
\address{Institute of Mathematics "Simion Stoilow",
Romanian Academy, P.O.Box 1-764, RO-70700,
Bucharest, Romania}
\email{Vasile.Brinzanescu@imar.ro} 

\thanks{The first author was partially supported by Swiss NSF contract SCOPES
2000-2003, No.7 IP 62615 and by contract CERES 39/2002-2004}

\author{Ruxandra Moraru}
\address{Department of Mathematics, University of Toronto, 100 St George Street, 
Toronto, Ontario, Canada, M5S 3G3}
\email{moraru@math.toronto.edu}


\begin{abstract}
In this paper, we study holomorphic rank-2 
vector bundles on non-K\" ahler elliptic surfaces. 
Our main tool for analysing these bundles is of course the spectral cover.
However, given the non-K\"{a}hler condition, the elliptic surfaces we are considering 
do not have sections and gerbes naturally arise in this context.
The spectral construction presented in this paper 
is a modification of the Fourier-Mukai transform for elliptic fibrations 
without a section. 
After examining some of the properties of this Fourier-Mukai transform,
we give a complete classification of vector bundles on these surfaces.
\end{abstract}

\maketitle

\section{Introduction}

The study of sheaves over elliptic fibrations has been a very active area of research
in both mathematics and physics over the past fifteen years; this is partly due to the role 
that such fibrations play in both mirror symmetry and the theory of integrable systems.
An object that has, more recently, proven very useful in their study is the Fourier-Mukai
transform. This transform is an equivalence of derived categories of sheaves
on elliptic fibrations with a section
whose properties are by now well-understood \cite{Mu,BBRP,B,BM};
for example, if $\mathcal{E}$ is a locally free sheaf on such a fibration $Y \rightarrow B$, 
then under some mild assumptions, the transform assigns
to $\mathcal{E}$ a torsion sheaf on the relative Jacobian $J(Y)$ of $Y$,
which is supported on its spectral cover.
An important point, which makes this construction possible, is the existence of a universal sheaf on 
$Y \times_B J(Y)$. 
However, if one considers an elliptic fibration without a section $X \rightarrow B$, 
then there is no universal sheaf on $X \times_B J(X)$; 
what exists instead is a twisted universal sheaf \cite{C1}, which can 
be used to define a ``twisted'' Fourier-Mukai transform 
that is now an equivalence of derived categories of twisted sheaves \cite{C2,DP}.

In this article, we consider holomorphic vector bundles on a specific class of 
fibrations without a section, non-K\" ahler elliptic surfaces.
Although bundles on projective elliptic fibrations have been extensively studied, 
not very much is known about the non-K\"{a}hler case;
another motivation for investigating bundles on these surfaces 
comes from recent developments in superstring theory,
where six-dimensional non-K\"{a}hler manifolds occur in the context of type IIA string
compactifications with non-vanishing background $H$-field
$-$ in fact, all the non-K\"{a}hler examples appearing in the physics literature so far
are non-K\"{a}hler principal elliptic fibrations (see \cite{BBDG,CCDLMZ,GP} and the references
therein).  
We first construct a particular twisted Fourier-Mukai transform
for locally free sheaves on non-K\"{a}hler principal elliptic bundles, 
transforming them into torsion sheaves, 
which has the advantage of allowing us to work with sheaves
instead of twisted sheaves; then, we use it to obtain a classification of
rank-2 vector bundles on arbitrary non-K\"{a}hler elliptic surfaces
(the existence and stability properties 
of such bundles are examined in \cite{Brinzanescu-Moraru1,Brinzanescu-Moraru2}). 
Note that this transform also makes sense for
coherent sheaves, but since we are primarily interested in classifying vector bundles,
we restrict our definition to locally free sheaves.
Furthermore, the techniques developed here naturally extend to 
the case of sheaves on higher dimensional 
non-K\"{a}hler elliptic and torus fibrations.  

A minimal non-K\"{a}hler elliptic surface $X$ is a Hopf-like surface that
admits a holomorphic fibration $\pi: X \rightarrow B$, over a smooth connected compact curve $B$, 
whose smooth fibres are isomorphic to a fixed smooth elliptic curve;
more precisely, if the surface $X$ does not have multiple fibres, 
then it is the quotient of a complex surface
by an infinite cyclic group.
The constructions presented here are based on methods used in \cite{Moraru}
to study bundles over Hopf surfaces.
If $X$ is a Hopf surface, then its Picard group is isomorphic to $\mathbb{C}^\ast$ and there 
exists a universal line bundle on $X \times \mathbb{C}^\ast$; one can then use this bundle
to define a natural transformation that takes a locally free sheaf $\mathcal{E}$ 
on $X$ to a torsion sheaf $\widetilde{\mathcal{L}}$ on $B \times \mathbb{C}^\ast$. Note that 
the relative Jacobian $J(X)$ of the Hopf surface $X$ is the quotient of $B \times \mathbb{C}^\ast$  
by the infinite cyclic group defining the surface;
however, the torsion sheaf $\widetilde{\mathcal{L}}$ cannot descend to this quotient. 
Nonetheless, this problem was solved in \cite{Moraru} by constructing a sheaf 
$\mathcal{N}$ on $B \times \mathbb{C}^\ast$ such that the tensor product
$\widetilde{\mathcal{L}} \otimes \mathcal{N}$ descends to $J(X)$. 
In the sequel, we show that this construction 
extends to any non-K\"{a}hler elliptic fibre bundle;
a twisted Fourier-Mukai transform is then defined as the composition of these two operations,
taking locally free sheaves on $X$ to torsion sheaves on $J(X)$. 

The article is organised as follows. 
We start by briefly reviewing the existence results that were proven in 
\cite{Brinzanescu-Moraru1}.
Then in the third section, we define the twisted Fourier-Mukai transform and 
examine some of its properties.
Next, we show that if a vector bundle $E$ is regular on all smooth fibres of $\pi$, then it is 
completely determined by its spectral cover $S_E$ and a certain line bundle on 
$S_E$ (over an elliptic curve, a bundle is said to be {\em regular} if 
its group of automorphisms is of the smallest possible dimension); 
consequently, we are able to prove that there is a one-to-one correspondence 
between rank-2 vector bundles with the same 
smooth spectral cover $S_E$ and a finite number of copies of a Prym variety 
associated to $S_E$. We end the article by giving an overview of the methods 
that can be used to classify bundles that are not regular on at least one
fibre of $\pi$.
\vspace{.1in}

{\bf Acknowledgements.} 
The first author would like to express his gratitude to the 
Max Planck Institute of Mathematics for its hospitality and 
stimulating atmosphere; part of this paper was prepared during his stay at
the Institute. The second author would like to thank Jacques Hurtubise for his generous 
encouragement and support during the completion of this paper.
She would also like to thank Ron Donagi and Tony Pantev for valuable discussions, and the
Department of Mathematics at the University of Pennsylvania for their hospitality,
during the preparation of part of this article.

\section{Holomorphic vector bundles}

Let $X\stackrel{\pi}{\rightarrow}B$ be a minimal non-K\"ahler elliptic 
surface, with $B$ a smooth compact connected curve; it is well-known that 
$X \stackrel{\pi}{\rightarrow} B$ is a quasi-bundle over $B$, that is, all the smooth fibres 
are pairwise isomorphic and the singular fibres are multiples of elliptic 
curves \cite{Kod,B2}.
Let $T$ be the general fibre of $\pi$, which is an elliptic curve, and
denote its dual $T^*$ (a non-canonical identification 
$T^*:= \mbox{Pic}^0(T)\cong T$). In this case, the Jacobian surface associated to 
$X\stackrel{\pi}{\rightarrow}B$  is simply 
\[ J(X)=B\times T^*\stackrel{p_1}{\rightarrow}B \]
(see, for example, \cite{Kod,BPV,B1}) and the surface 
$X$ is obtained from its relative Jacobian by a finite number 
of logarithmic transformations \cite{Kod,BPV,BrU}.

\subsection{Line bundles} 
\label{line bundles}
Before giving a general description of line bundles on $X$, we begin with torsion line bundles.
Suppose that $\pi$ has a multiple fibre $mF$ over the point $b$ in $B$; the line bundle 
associated to the divisor $F$ of $X$ is then such that 
$(\mathcal{O}_X(F))^m = \mathcal{O}_X(mF) = \pi^\ast\mathcal{O}_B(b)$.
Let $P_2$ denote the subgroup of $\text{Pic}(X)$ generated by $\pi^\ast \text{Pic}(B)$ and
the $\mathcal{O}_X(T_i)$', where $T_1, \dots, T_l$ are the multiple fibres (if any) of $X$. 
The group of all torsion line bundles on  $X$ is then given by
\begin{equation}\label{torsion line bundles} 
\text{Pic}^\tau(X) = P_2 \otimes \mathbb{C}^\ast. 
\end{equation}
If $X$ does not have multiple fibres, then
the set of all holomorphic line bundles on $X$ with 
trivial Chern class is given by the zero component of the Picard group
$\text{Pic}^0(X) \cong \text{Pic}^0(B) \times \mathbb{C}^\ast$.
In this case, any line bundle in $\text{Pic}^0(X)$ is of the form
$H \otimes L_\alpha$,
where $H$ is the pullback to $X$ of an element of $\text{Pic}^0(B)$ and
$L_\alpha$ is the line bundle corresponding to the
constant automorphy factor $\alpha \in \mathbb{C}^\ast$;
in particular, there exists a universal 
(Poincar\'{e}) line bundle $\mathcal{U}$ on $X \times \text{Pic}^0(X)$ whose restriction to
$X \times \mathbb{C}^\ast := X \times \{ 0 \} \times \mathbb{C}^\ast$ is constructed in terms of
constant automorphy factors (for details, see \cite{Brinzanescu-Moraru1}). 

For vector bundles on any elliptic surface $X$, restriction to a fibre is a natural operation; 
however, it is important to note that if $X$ is non-K\"{a}hler,
then the restriction of {\em any} line bundle on $X$ to a smooth fibre of $\pi$ {\em always} 
has degree zero \cite{Brinzanescu-Moraru1}.
Furthermore, even though non-K\"{a}hler elliptic surfaces have very few divisors
(they are given by the fibres of $\pi$), 
there exist many line bundles on them;
we have the following classification \cite{Brinzanescu-Moraru1}. 

\begin{proposition}
Let $X \stackrel{\pi}{\rightarrow} B$ be a non-K\"{a}hler elliptic surface with 
general fibre $T$ and $J(X)=B\times T^*$ be its relative Jacobian. 
Fix a section $\Sigma \subset J(X)$. Then:

(i) There exists a line bundle 
$L$ on $X$ whose restriction to every smooth fibre $T_b = \pi^{-1}(b)$ of $\pi$ is the same 
as the line bundle $\Sigma_b$ of degree zero on $T \cong T_b$.

(ii) The set of all line bundles on $X$ whose restriction to every smooth fibre of $\pi$
is determined by the section $\Sigma$ is a principal homogeneous space over $P_2$.
\end{proposition}

\subsection{The spectral construction}
\label{spectral curve}
Consider a pair $(c_1,c_2) \in NS(X) \times \mathbb{Z}$; then, its 
{\em discriminant} is defined as
\[ \Delta(2,c_1,c_2):=\frac{1}{2} \left( c_2 - \frac{c_1^2}{4} \right).\]
Let $E$ be a rank 2 vector bundle over $X$, with $c_1(E) = c_1$ and 
$c_2(E) = c_2$. 
For the remainder of the paper, we fix the notation:
\[ \Delta(E) := \Delta(2,c_1,c_2). \]
To study bundles on $X$, one of our main tools is 
restriction to the smooth fibres of the fibration $\pi: X \rightarrow B$. 
Since the restriction of any bundle on $X$ 
to a fibre $T$ has first Chern class zero, we can consider $E$ as a 
family of degree zero bundles
over the elliptic curve $T$, parametrised by $B$. These families were described in detail 
in \cite{Brinzanescu-Moraru1}. For the convenience of the reader, we give a brief summary.

Given a rank two bundle over $X$, its restriction to a generic
fibre of $\pi$ is semistable. More precisely, its restriction to a fibre
$\pi^{-1}(b)$ is unstable on at most an isolated set of points $b \in B$;
these isolated points are called the {\em jumps} of the bundle.
Furthermore, there exists a divisor in the relative Jacobian 
$J(X) = B \times T^\ast$ of $X$, 
called the {\em spectral curve} or {\em cover} of the bundle, 
that encodes the isomorphism
class of the bundle over each fibre of $\pi$. 
The spectral curve can be constructed as follows.
If the surface $X$ does not have multiple fibres, then
there exists a universal bundle $\mathcal{U}$ on $X \times \text{Pic}^0(X)$,
whose restriction to $X \times \mathbb{C}^\ast$ is also denoted $\mathcal{U}$;
we associate to the rank-2 vector bundle $E$ the sheaf on $B \times \mathbb{C}^\ast$ defined by 
\[ \widetilde{\mathcal{L}} := R^1\pi_\ast (s^\ast E \otimes  \mathcal{U}), \] 
where $s:X \times \mathbb{C}^\ast \rightarrow X$ 
is the projection onto the first factor, $id$ is the identity map, and
$\pi$ also denotes the projection 
$\pi := \pi \times id : X \times \mathbb{C}^\ast \rightarrow B \times \mathbb{C}^\ast$. 
This sheaf is supported on 
a divisor $\widetilde{S_E}$, defined with multiplicity, 
that descends to a divisor $S_E$ in $J(X)$ of the form
\[ S_E := \left( \sum_{i=1}^k \{ x_i \} \times T^\ast \right) + \overline{C}, \]
where $\overline{C}$ is a bisection of $J(X)$ and $x_1, \cdots, x_k$ are points in $B$
that correspond to the jumps of $E$. 
The spectral curve of $E$ is defined to be the divisor $S_E$.

If the fibration $\pi$ has multiple fibres, then one can associate to 
$X$ a principal $T$-bundle $\pi': X' \rightarrow B'$ over an $m$-cyclic covering 
$\varepsilon: B' \rightarrow B$, where the integer $m$ depends on the multiplicities
of the singular fibres;
note that the map $\varepsilon$ induces natural $m$-cyclic coverings
$J(X') \rightarrow J(X)$ and $\psi :X'\rightarrow X$. 
By replacing $X$ with $X'$ (which does not have multiple fibres)  
in the above construction, we obtain the  spectral cover $S_{\psi^\ast E}$ of the 
vector bundle $\psi^\ast E$ as a divisor in $J(X')$.
We then define the spectral cover $S_E$ of $E$ to be the projection of $S_{\psi^\ast E}$ in $J(X)$.

\begin{note}
The spectral construction can be defined for any rank $r$ vector bundle. In particular, for a line
bundle, the spectral cover corresponds to the section of the Jacobian surface $J(X)$ defined in 
section \ref{line bundles}.
\end{note}

\begin{remark}\label{graph}
Let $\delta$ be a line bundle on $X$. It then defines the an involution $i_\delta$ of
the relative Jacobian $J(X) = B \times T^\ast$ of $X$, given by
$(b,\lambda) \mapsto (b,\delta_b \otimes \lambda^{-1})$, 
where $\delta_b$ denotes the restriction of $\delta$ to the fibre 
$T_b = \pi^{-1}(b)$.  
Taking the quotient of $J(X)$ by this involution, we obtain a ruled surface 
$\mathbb{F}_\delta := J(X)/i_\delta$ over $B$. 
Let $\eta: J(X) \rightarrow \mathbb{F}_\delta$ be the canonical map.
By construction, the spectral curve $S_E$ of a bundle $E$ 
with determinant $\delta$ is invariant with respect to
the involution $i_\delta$ and descends to the quotient 
$\mathbb{F}_\delta$; in fact, it is the pullback via $\eta$ of a 
divisor on $\mathbb{F}_\delta$
of the form $\mathcal{G}_E := \sum_{i=1}^k f_i + A$,
where $f_i$ is the fibre of the ruled surface $\mathbb{F}_\delta$ over 
the point $x_i$ and
$A$ is a section of the ruling such that $\eta^\ast A = \overline{C}$. 
\end{remark}

We have seen that one can associate to every vector bundle on $X$ a spectral cover.
Conversely, given a fixed spectral cover without vertical components, 
there exists at least one rank-2 vector bundle on $X$ associated to it. 
More precisely, we have the following existence result
\cite{Brinzanescu-Moraru1}.

\begin{theorem}
\label{existence}
Let $X$ be a non-K\"{a}hler elliptic surface and fix a line bundle $\delta$ on $X$. 
Moreover, let $i_\delta$ be the involution of the relative Jacobian $J(X)$ determined by $\delta$ 
and suppose that $\overline C$ is a bisection of $J(X)$ that is invariant with
respect to $i_\delta$. Then, there exists at least one holomorphic rank-2 vector bundle on $X$ 
that has determinant $\delta$ and spectral cover $\overline{C}$.
\end{theorem}

\begin{note}
Further existence results are given in \cite{Brinzanescu-Moraru1}.
\end{note}

\section{Gerbes and the Fourier-Mukai transform}
\label{gerbe}
For elliptic fibrations with a section, the spectral construction  
has, more recently, been studied from the point of view of the Fourier-Mukai transform.
Let us consider for a moment an elliptic fibration $Y \stackrel{\pi}{\rightarrow} B$ that has
a section; in this case, the section of the fibration gives a natural identification
of $Y$ with its relative Jacobian $J(Y)$, that is, there is a canonical isomorphism $J(Y) \cong Y$.
Furthermore, there exists a universal line bundle $\mathcal{P}$ on the fibred product 
$Y \times_B J(Y) = Y \times_B Y$.
Denote $p_1$ and $p_2$ the projections of $Y \times_B Y$ onto the first and second factors, 
respectively, and let $E$ be a holomorphic vector bundle on $Y$. As in the previous section, the spectral 
cover $C$ of $E$ can be described as the support in 
$Y \cong J(Y)$ of the torsion sheaf $R^1{p_2}_\ast(p^\ast_1E \otimes \mathcal{P})$,
which is the Fourier-Mukai transform of $E$, under some mild assumptions on $E$.
Moreover, if $E$ is semistable of degree zero on the generic fibre of $\pi$, then the
restriction of $R^1{p_2}_\ast(p^\ast_1E \otimes \mathcal{P})$ to $C$ is a rank one coherent
torsion free sheaf $L$ and one can recover the bundle $E$ from the spectral data $( C,L )$.

Let us turn to the case of a non-K\"{a}hler elliptic fibre bundle $X \stackrel{\pi}{\rightarrow} B$;
given the non-K\"{a}hler condition, $\pi$ is now an elliptic fibration 
without a section $-$ a genus one fibration in the terminology of \cite{DP}.
Consequently, $X$ is not isomorphic to its
relative Jacobian; indeed, if the general fibre of $X$ is the elliptic curve $T$,
then $J(X) = B \times T^\ast$. Furthermore,
there is no universal line bundle on $X \times_B J(X)$. 
To construct spectral covers of vector bundles on $X$, 
we bypassed this problem by working instead with the universal line bundle
$\mathcal{U}$ that exists on $X \times \mathbb{C}^\ast$
(see sections \ref{line bundles} and \ref{spectral curve}).
Consider a rank-2 vector bundle $E$ on $X$.
In this section, we show that the spectral cover $S_E \subset J(X)$ of $E$ can also be 
determined as the support of a torsion sheaf 
$\mathcal{L}$ on $J(X)$, leading us to the natural definition of a twisted Fourier-Mukai 
transform $\Phi$ for locally free sheaves on $X$; in particular, $\Phi(E) = \mathcal{L}$.
Some properties of this transform are then examined.

\subsection{Gerbes}
\label{gerbes}

Let $X$ be a general genus one fibration without multiple fibres.
Then, as in the case of non-K\"{a}hler elliptic surfaces, $X$ is not isomorphic
to its relative Jacobian $J(X)$ and distinct fibrations may have the same relative Jacobian;
the information lost
by passing from $X$ to its Jacobian is, however, encoded in 
the class of $X$ in the Tate-Shavarevich group of $J(X)$, or equivalently,
by an $\mathcal{O}_{J(X)}^\ast$-gerbe $\xi$ on $J(X)$. 
Furthermore, even though there is no universal sheaf on $X \times_B J(X)$, 
universal $\xi$-twisted sheaves do exist \cite{C1}.
An equivalence between derived categories of twisted sheaves,
called a ``gerby'' Fourier-Mukai transform, 
can then be constructed by using these universal twisted sheaves
(see \cite{DP} for precise statements).

For example, if $X$ is a non-K\"{a}hler
principal $T$-bundle over the curve $B$, 
then the $\mathcal{O}_{J(X)}^\ast$-gerbe corresponding to $X$ is given by a 
cohomology class $\xi \in H^1(B,\mathcal{O}_B(T))$, where $\mathcal{O}_B(T)$ denotes the
sheaf of germs of locally holomorphic maps from $B$ to $T$, 
and there exists a $\xi$-twisted
universal sheaf $\mathcal{U}_\xi$ on $X \times_B J(X) = X \times T^\ast$. 
In this case, the results of \cite{DP} state that the twisted Fourier-Mukai transform
defined by $\mathcal{U}_\xi$ is an equivalence between the derived category $D^b(J(X),\xi)$ 
of $\xi$-twisted sheaves on $J(X)$ and the derived category $D^b(X)$ of sheaves on $X$.
But in this paper, we use the gerbe $\xi$ to construct
an explicit twist of the torsion sheaf  
$\widetilde{\mathcal{L}} := R^1\pi_\ast (s^\ast E \otimes  \mathcal{U})$, on $B \times \mathbb{C}^\ast$,
that descends to $J(X) = B \times T^\ast$.

The cohomology class $\xi$ induces an isomorphism of the surface $X$
with a quotient of the form
\[ X = \Theta^\ast / \langle \tau \rangle ,\]
where $\Theta$ is a line bundle on $B$ with positive Chern class $d$,
$\Theta^\ast$ is the complement of the zero section in the total space of 
$\Theta$, and $\langle \tau \rangle$ is the multiplicative cyclic group 
generated by a fixed complex number $\tau \in \mathbb{C}$, with 
$| \tau | > 1$; given this identification, every fibre of $\pi$ is isomorphic to the elliptic
curve $\mathbb{C}^\ast/\tau^n$. 
Hence, multiplication by $\tau$ defines a natural $\mathbb{Z}$-action on 
$X \times \mathbb{C}^\ast$ that is trivial on $X$, inducing the quotient
$(X \times \mathbb{C}^\ast)/\mathbb{Z} = X \times T^\ast = X \times_B J(X)$.

We saw in section \ref{line bundles} that every torsion line bundle
$L \in \text{Pic}^0(X)$ decomposes uniquely 
as $L = H \otimes L_\alpha$, for $H \in \pi^\ast \text{Pic}^0(B)$
and $\alpha \in \mathbb{C}^\ast$. 
However, if one considers the constant factor of automorphy $\tau \in \mathbb{C}^\ast$, 
then one easily verifies that $L_\tau$ is in fact the pullback of a line bundle on $B$ of degree $-d$,
also denoted $L_\tau$. 
Note that there is also a natural $\mathbb{Z}$-action on $B \times \mathbb{C}^\ast$,
defined as multiplication by $\tau$ on the second factor, and 
$ (B \times \mathbb{C}^\ast)/\mathbb{Z} \cong J(X)$. 
Moreover, this action extends to the torsion sheaf 
$\widetilde{\mathcal{L}} := R^1\pi_\ast (s^\ast E \otimes  \mathcal{U})$,
taking the stalk $\widetilde{\mathcal{L}}_{(x,\alpha)}$ to
$\widetilde{\mathcal{L}}_{(x,\tau\alpha)} \otimes L_{\tau^{-1},x}$ \cite{Brinzanescu-Moraru1}.
Therefore, $\widetilde{\mathcal{L}}$ cannot  descend to $J(X)$ because 
it is not invariant with respect to this action.
To fix this problem, we construct a sheaf $\mathcal{N}$ on  
$B \times \mathbb{C}^\ast$ and a $\mathbb{Z}$-action that leaves
the tensor product $\widetilde{\mathcal{L}} \otimes \mathcal{N}$ invariant. 

Let us choose a point $b_0$ in $B$ over which the graph of $E$ does not have a 
vertical component over; without loss of generality, 
we assume that the line bundle $L_{\tau}$ is given by the divisor $-db_0$
(if the divisor associated to $L_\tau$ is instead $\sum_i b_i - \sum_j b_j$, the construction 
below follows verbatim, but the notation is more complicated). 
Let $W = (b_0\times \mathbb{C}^\ast) \cap \widetilde{S_E}$ be the set of points
on $\widetilde{S_E}$ lying above $b_0$. 
If $(a,b)$ is a representation of the pair of points on $S_E$ above $b_0$, 
then $W$ is the set of all  translates of this pair by 
$\tau$, that is, $W = \bigcup_{i \in \mathbb{Z}} (\tau^ia , \tau^ib)$.
Also, $(a+b)$ is a divisor on $\widetilde{S_E}$ and
we denote $T^i(a+b) := \tau^ia + \tau^ib$ the translate of $(a+b)$ 
by $\tau^i$; we define a divisor on $\widetilde{S_E}$ as the locally 
finite sum
\[ D_{E,b_0} := \sum_{i \in \mathbb{Z}} idT^i(a+b).\]
Let $\mathcal{N}$ be the line bundle on $\widetilde{S_E}$ associated to the 
invertible sheaf $\mathcal{O}_{\widetilde{S_E}}(D_{E,b_0})$; we also denote by 
$\mathcal{N}$ the line bundle considered as a sheaf on 
$B \times \mathbb{C}^\ast$. 

Fix a section $\gamma$ of the line bundle $L_{\tau}$;
we use it to define the following $\mathbb{Z}$-action 
on the sheaf
$\widetilde{\mathcal{L}} \otimes \mathcal{N}$ over $B \times \mathbb{C}^\ast$:
\begin{diagram}
\widetilde{\mathcal{L}} \otimes \mathcal{N}|_{(z,\alpha)} & 
\rTo_{\tau}^{(z,\alpha,s \otimes t)  \mapsto  
(z,\tau \alpha, \alpha s \otimes \gamma t)} & 
\widetilde{\mathcal{L}} \otimes \mathcal{N}|_{(z,\tau \alpha)}\\
\dTo &  & \dTo\\
B \times \mathbb{C}^\ast  & \rTo^{\tau}_{(z,\alpha)  \mapsto 
(z,\tau \alpha)} &  B \times \mathbb{C}^\ast.
\end{diagram}
Clearly, the sheaf $\widetilde{\mathcal{L}} \otimes \mathcal{N}$ is invariant under this action 
and descends to the quotient
$J(X)$;
we denote the quotient sheaf
\[ \mathcal{L} := (\widetilde{\mathcal{L}} \otimes \mathcal{N})/\sim.\]
Note that the support of $\mathcal{L}$ is $S_E$;
moreover, 
if we take the pull back of $\mathcal{L}$ to 
$B \times \mathbb{C}^\ast$ and tensor it by $\mathcal{N}^\ast$, then
we recover $\widetilde{\mathcal{L}}$ (as above, we also denote $\mathcal{N}^\ast$ 
the sheaf on $B \times \mathbb{C}^\ast$ obtained by extending the line bundle $\mathcal{N}^\ast$ 
on $\widetilde{S_E}$ by zero outside $\widetilde{S_E}$).

\subsection{Properties of the Fourier-Mukai transform}

Consider the following commutative diagram:
\begin{diagram}[size=2em]
 X \times_B J(X) & \lTo^{q} & X \times \mathbb{C}^\ast & \rTo{\pi} & B \times \mathbb{C}^\ast
& \rTo^{\rho} & J(X)\\ 
 \dTo<{\underline{s}} &  & \dTo<{s} & & \dTo & & \dTo \\
X & \rEqual & X & \rTo_{\pi} & B & \rEqual & B,
\end{diagram}
where $\underline{s}$ is projection onto the first factor,
$q : X \times \mathbb{C}^\ast \rightarrow X \times T^\ast = X \times_B J(X)$
and $\rho: B \times \mathbb{C}^\ast \rightarrow B \times T^\ast = J(X)$ 
are the natural quotient maps induced by the $\mathbb{Z}$-actions 
defined in section \ref{gerbes}, 
and $\pi$ and $s$ are the projections defined in section \ref{spectral curve}.
Given a locally free sheaf $E$ on $X$, we define the twisted Fourier-Mukai transform 
to be the complex of sheaves $\Phi(E)$ on $J(X)$ given by
\[\Phi(E) := 
\left( R\pi_\ast \left(s^\ast E \otimes \mathcal{U} \right) \otimes \mathcal{N} \right)/\sim. \]
Conversely, if $\mathcal{L}$ is a sheaf on $J(X)$, we define the ``inverse'' twisted Fourier-Mukai
transform as the complex of sheaves $\hat{\Phi}(\mathcal{L})$ on $X$ given by
\[ \hat{\Phi}(\mathcal{L}) := 
R\underline{s}_\ast \left( \left( \pi^\ast \left( \left( \rho^\ast \mathcal{L} \right) 
\otimes \mathcal{N}^\ast \right) \otimes \mathcal{U}^\ast \right) / \sim \right). \]
Referring to section \ref{gerbes}, the sheaves of these complexes are well-defined;
however, we should point out that we use the term ``inverse'' only in the sense that
the transforms $\Phi$ and $\hat{\Phi}$ are inverses of each other on certain classes of sheaves, 
as will become clear in the following.
We state some of their properties in:
\begin{theorem}\label{fourier-mukai}
(i) Suppose that $E$ is a rank-2 vector bundle on $X$ without jumps. Then,
$\Phi^0(E) = 0$ and $\hat{\Phi}^0(\Phi^1(E)) = E$.

(ii) If $\mathcal{L}$ is a torsion sheaf on $J(X)$, supported on a bisection $C \subset J(X)$, 
that has rank 1 on the smooth points of $C$ and rank at most 2 on the singular ones,
then $\hat{\Phi}^1(\mathcal{L}) = 0$ and $\Phi^1(\hat{\Phi}^0(\mathcal{L})) = \mathcal{L}$.
\end{theorem}

\begin{proof}
Let $E$ be a rank-2 vector bundle on $X$ without jumps.
If one unravels the definitions, then one sees that
\[ \hat{\Phi}^0(\Phi^1(E)) = \underline{s}_\ast \left( \left( \left( \pi^\ast R^1\pi_\ast \left(
s^\ast E \otimes \mathcal{U} \right) \right) \otimes \mathcal{U}^\ast \right) / \sim \right); \]
moreover, flat base change induces a natural isomorphism 
$\pi^\ast R^1\pi_\ast(\mathcal{F}) = \mathcal{F}$ for any coherent sheaf $\mathcal{F}$ on 
$B \times \mathbb{C}^\ast$, implying that  
$\hat{\Phi}^0(\Phi^1(E)) = \underline{s}_\ast(s^\ast E / \sim )$.
However, since $E$ is invariant under the $\mathbb{Z}$-action, we have 
$s^\ast E / \sim  = \underline{s}^\ast E$
and $\hat{\Phi}^0(\Phi^1(E)) = \underline{s}_\ast(\underline{s}^\ast E) \cong E$,
proving (i).

Note that for any coherent sheaf $\mathcal{G}$ on $X \times_B J(X)$,
there is an identification 
$s^\ast \underline{s}_\ast(\mathcal{G}) = q^\ast(\mathcal{G})$,
which is again due to flat base change.
Consequently, if $\mathcal{L}$ is a torsion sheaf on $J(X)$ satisfying  
the hypothesis of part (ii),
we have 
\[ \Phi^1(\hat{\Phi}^0(\mathcal{L})) = 
\left( R^1\pi_\ast \left( \pi^\ast \left( \rho^\ast \mathcal{L} \otimes \mathcal{N}^\ast 
\right) \right) \otimes \mathcal{N} \right) / \sim . \]
There is a natural map 
$\rho^\ast \mathcal{L} \otimes \mathcal{N} \rightarrow 
R^1\pi_\ast(\pi^\ast(\rho^\ast \mathcal{L} \otimes \mathcal{N}))$;
since both sheaves are supported on $\rho^\ast C$ and have naturally isomorphic fibres,
this map must be an isomorphism, giving
$\Phi^1(\hat{\Phi}^0(\mathcal{L})) = 
\left( \left(\rho^\ast \mathcal{L} \otimes \mathcal{N}^\ast 
\right) \otimes \mathcal{N} \right) / \sim = \mathcal{L}$.
\end{proof}

\begin{remark}\label{uniqueness of line bundle}
Suppose that the bisection $C$ is smooth. 
If the torsion sheaf $\mathcal{L}$ satisfies the conditions of Theorem \ref{fourier-mukai} (ii), 
then its restriction to $C$ is a line bundle.
Denote $L := \left( \pi^\ast \left( \left( \rho^\ast \mathcal{L} \right) 
\otimes \mathcal{N}^\ast \right) \otimes \mathcal{U}^\ast \right) / \sim$
the quotient sheaf on $X \times_B J(X)$; 
one easily verifies that the support of $L$ is the principal $T$-bundle 
$W := X \times_B C$ and that the restriction of $L$ to $W$ is a line bundle.
Let $\tilde{\gamma} : W \rightarrow X$ be the natural projection. 
We then have the following identification:
\[ \hat{\Phi}^0(\mathcal{L}) = \underline{s}_\ast (L) = \tilde{\gamma}_\ast(L|_W).\]
In particular, we see that if $E$ is a rank-2 vector bundle with spectral cover $C$
and $\Phi^1(E) = \mathcal{L}$, then $L_W := L|_W$ is the unique line bundle on $W$ such that
$\tilde{\gamma}_\ast(L_W) = E$.
\end{remark}

\section{Regular rank two vector bundles}
\label{regular rank-2 bundles} 

\subsection{Definition and existence}

Over a smooth elliptic curve, a rank-2 vector bundle of degree zero 
is said to be {\em regular} if its group of automorphisms is of the smallest dimension,
that is, it is a semistale bundle that is
{\em never} isomorphic to $\lambda_0 \oplus \lambda_0$, with 
$\lambda_0 \in \text{Pic}^0(T)$.
If $E$ is a holomorphic rank-2 vector bundle on an elliptic surface $X$, then $E$
is regular if and only if its restriction to a fibre 
$T$ is always regular. 
One can easily show that if a rank-2 vector bundle has a smooth spectral cover,
then it is regular (see, for example, \cite{FMW,Moraru});
given Theorem \ref{existence}, regular bundles therefore always exits.
However, regular bundles do not always have smooth spectral covers, as stated in:
\begin{proposition}\label{existence of regular bundles}
Let $X$ be a non-K\"{a}hler elliptic fibre bundle over the curve $B$ 
and let $\delta$ be a line bundle on $X$. 
Consider a (not necessarily smooth) bisection $C$ of the relative Jacobian $J(X)$ of $X$ that is 
invariant with respect to the involution $i_\delta$ of $J(X)$ determined by $\delta$.
Then, there exists at least one {\em regular} bundle on $X$ with determinant $\delta$ and 
spectral cover $C$.
\end{proposition}
\begin{proof}
Let us consider the section $A$ of 
$\mathbb{F}_\delta$ corresponding to $C$ (see remark \ref{graph}).
Locally, over an open disc $D$ of $B$, this section is the graph of a rational map
$\varrho: D \rightarrow \mathbb{P}^1$. Choose a line bundle $\lambda$ over $T$ 
with $c_1(\lambda) = -1$; non-trivial extensions of $\lambda^\ast$ by $\lambda$ 
are therefore regular rank-2 vector bundles over $T$, which are parametrised 
by $\mathbb{P}(H^1(T,\lambda^2)) \cong \mathbb{P}^1$. 
If we write $\varrho$ as a quotient $p(b)/q(b)$, where $p$ and $q$ are
coprime polynomials, we can define a regular bundle on $D \times T$, whose graph is the 
restriction of $A$ to $D$, by the extension
\begin{equation}\label{regular extension} 
\begin{array}{rcl}
D & \longrightarrow & H^1(T,\lambda^2),\\
b & \longmapsto & (p(b),q(b)). \end{array} 
\end{equation}
Moreover, if two regular bundles on $D \times T$ have the same graph,
given by a rational map $\varrho: D \rightarrow \mathbb{P}^1$, then they are
isomorphic (see Lemma 5.1.2 of \cite{B-H}); this is also true for any simply connected 
subset $D$ of $B$, in which case the isomorphism can be chosen to be
an $SL(2,\mathbb{C})$-isomorphism. Let $\{ D_i \}$ be an open cover of $B$ such that
every open set $D_i$ is a disc in $B$ and the intersection of any two open sets is 
simply connected. The extension \eqref{regular extension} gives regular
bundles over each $D_i$, which can then be glued together to give a regular bundle
on $X$ with the required properties.
\end{proof}

\subsection{Classification of regular bundles}
Assume that the surface $X$ has multiple fibres $m_1T_1, \dots, m_lT_l$ over the
points $b_1, \dots, b_l$ in $B$. We fix a line bundle $\delta$ on $X$ and its associated
involution $i := i_\delta$ of the Jacobian surface $J(X)$.
Consider a smooth irreducible bisection $C$ of $J(X)$, 
invariant under the involution $i$, and let $W$ be the normalisation of $X \times_B C$.
As $C$ is a double cover of $B$, then $W$ is also a double cover of $X$ and 
we have the following commutative diagram
\begin{diagram}[size=2em]
 W & \rTo^{\tilde{\gamma}} & X \\
\dTo<{\rho} & & \dTo>{\pi} \\
C & \rTo_{\gamma} & B.
\end{diagram}
Note that, to each multiple fibre $m_iT_i$ of $\pi$, there corresponds a pair of multiple fibres 
$m_iT'_i$ and $m_iT''_i$ of $\rho$ (that both get mapped to $m_iT_i$ under $\tilde{\gamma}$);
furthermore, if the fibre $m_jT_j$ lies over a branch point of $C \stackrel{\gamma}{\rightarrow} B$,
then $m_jT'_j = m_jT''_j$.
\newline
\indent
Let $P_{2,W}$ be the subgroup of $\rm{Pic}(W)$ generated by $\rho^\ast \rm{Pic}(C)$, 
$\mathcal{O}_W(T'_i)$, and  $\mathcal{O}_W(T''_i)$, $i = 1, \dots, l$; we define the following 
subgroup of $P_{2,W}$:
\begin{equation}\label{subgroup1}
P_{2,W}^0 := \left\{ 
\mbox{$\lambda \in P_{2,W}$ : $\tilde{\gamma}_\ast(c_1(\lambda)) = 0$ in $H^2(X,\mathbb{Z})$} \right\}.
\end{equation}
Finally, suppose that there are $s$ multiple fibres $m_{i_1}T_{i_1}, \dots, m_{i_s}T_{i_s}$
of $\pi$ that do not lie over branch points of $C \stackrel{\gamma}{\rightarrow} B$. 
Then, we denote 
\begin{equation}\label{subgroup2}
P_{2,W}^{T' - T''} := \langle \mathcal{O}_W(T'_{i_j} - T''_{i_j}), j=1, \dots, s \rangle \subset P_{2,W}
\end{equation}
the finite subgroup of order $m_{i_1} \dots m_{i_s}$
generated by line bundles on $W$ of the form $\mathcal{O}_W(T'_{i_j} - T''_{i_j})$.
Given the above notation, we can now state the main result of this section.

\begin{theorem}\label{classification of regular bundles}
Let us denote $\Sigma_C$ the section of $J(W) = C \times T^\ast$ that corresponds to the map 
$C \rightarrow J(W)$. Fix a line bundle $L$ on $W$ that corresponds to the section $\Sigma_C$
and is such that $\det(\tilde{\gamma}_\ast(L)) = \delta$.
All rank-2 vector bundles on $X$ with determinant $\delta$ and 
spectral cover $C$ are then given precisely by
\[ \tilde{\gamma}_\ast(L \otimes \lambda), \ \ \lambda \in P_{2,W}^0/\bar{\iota}, \]
where $\bar{\iota}$ denotes the involution on $W$ that interchanges the sheets of $\tilde{\gamma}$
and 
\[ P_{2,W}^0/\bar{\iota} := \left\{ \lambda \in P_{2,W}^0 \ | \ 
\bar{\iota}^\ast \lambda \otimes \lambda = \mathcal{O}_W \right\}.\]
Moreover, there is a natural exact sequence
\[ 0 \rightarrow Prym(C/B) \rightarrow P_{2,W}^0/\bar{\iota} \rightarrow P_{2,W}^{T' - T''} 
\rightarrow 0,\]
where $Prym(C/B)$ is the Prym variety associated to the double cover 
$C \stackrel{\gamma}{\rightarrow} B$.
Note that if $X$ does not have multiple fibres, then $P_{2,W}^0/\bar{\iota} \cong Prym(C/B)$.
\end{theorem}

\begin{remark*}
If one does not fix the determinant, rank-2 vector bundles on $X$ 
with first Chern class $c_1(\delta)$ and 
spectral cover $C$ are then parametrised by the group $P_{2,W}^0$,
which is isomorphic to  $m_{i_1}m_{i_2} \dots m_{i_s}$ copies of the Jacobian $J(C)$ of $C$.
\end{remark*}

\begin{proof}
Consider a rank-2 vector bundles $E$ with determinant $\delta$
and spectral cover $C$. Referring to remark \ref{uniqueness of line bundle}, 
if $X$ does not have multiple
fibres, we can associate to $E$ a unique line bundle $L$ on $W$ such that $E = \tilde{\gamma}_\ast(L)$.
If $X$ has multiple fibres, 
the same can be said of the pullback $\psi^\ast E$ of $E$ to $X'$,
where $X'$ is the $m$-to-one cover of $X$ (without multiple fibres)
described in section \ref{spectral curve}.
Indeed, if $N_1$ and $N_2$ are non-isomorphic line bundles on $W$, they must generate 
non-isomorphic rank-2 vector bundles on $X$, otherwise, their pullbacks to $X'$ 
would contradict the above.
In other words, whether $X$ has multiple fibres or not,
there is a one-to-one correspondence between rank-2 vector bundles, with 
determinant $\delta$ and spectral cover $C$, and 
line bundles $N$ on $W$, associated to the section $\Sigma_C \subset J(W)$,
such that $\det(\tilde{\gamma}_\ast(N)) = \delta$.

Fix a line bundle $L$ on $W$ that corresponds to the section $\Sigma_C$
and is such that $\det(\tilde{\gamma}_\ast(L)) = \delta$.
Recall that any line bundle corresponding to $\Sigma_C$ is of the form $L \otimes \lambda$,
where $\lambda$ is an element of $P_{2,W}$.
Hence, consider a line bundle $L' = L \otimes \lambda$, $\lambda \in P_{2,W}$, 
such that $\det(\tilde{\gamma}_\ast(L')) = \delta$;
in particular, $c_1(\tilde{\gamma}_\ast(L')) = c_1(\delta)$.
Since 
\[ c_1(\tilde{\gamma}_\ast(L \otimes \lambda)) = c_1(\tilde{\gamma}_\ast(L))
+ \tilde{\gamma}_\ast(c_1(\lambda)),\]
we see that $c_1(\tilde{\gamma}_\ast(L')) = c_1(\delta)$ if and only if 
$\tilde{\gamma}_\ast(c_1(\lambda)) = 0$, implying that $\lambda \in P_{2,W}^0$.
For any line bundle $N$ on $W$, there is an exact sequence on $W$
\begin{equation}\label{exact}
0 \rightarrow \bar{\iota}^\ast N \otimes \tilde{\gamma}^\ast L_0^{-1} \rightarrow 
\tilde{\gamma}^\ast \tilde{\gamma}_\ast(N) \rightarrow N \rightarrow 0.
\end{equation}
Inserting the line bundles $L$ and $L'$ in the exact sequence \eqref{exact}, 
we have
\[ \det(\tilde{\gamma}^\ast\tilde{\gamma}_\ast(L')) =
\det(\tilde{\gamma}^\ast\tilde{\gamma}_\ast(L)) \otimes (\bar{\iota}^\ast \lambda \otimes \lambda). \]
Consequently, $\lambda$ must satisfy the equation
\begin{equation}\label{involution invariance}
\bar{\iota}^\ast \lambda \otimes \lambda = \mathcal{O}_W.
\end{equation}
The set of all rank-2 vector bundles on $X$ with determinant $\delta$ and spectral cover $C$
is therefore parametrised by $P_{2,W}^0/\bar{\iota}$.

We end by determining the generators of the group $P_{2,W}^0/\bar{\iota}$.
Let $\lambda \in P_{2,W}^0$.
One easily verifies that if $\lambda \in \rho^\ast \rm{Pic}(C)$, then 
$\tilde{\gamma}_\ast(c_1(\lambda)) = 0$ if and only if $\lambda \in \rho^\ast \rm{Pic}^0(C)$. 
We then assume that the divisor of the line bundle $\lambda$
contains multiple fibres of $W \stackrel{\rho}{\rightarrow} C$. 
Clearly, if it contains the multiple $a_iT'_i$, with $0 \leq a_i \leq m_i-1$,
then $\tilde{\gamma}_\ast(c_1(\lambda)) = 0$ if and only if it also contains $a_iT''_i$.
Hence, the divisor of $\lambda$ must contain multiples of the form 
$a_1(T'_1 - T''_1) + \dots + a_r(T'_l - T''_l)$, with $0 \leq a_i \leq m_i - 1$, 
for all $i = 1, \dots, l$.
Since $T'_j = T''_j$ whenever $T_j$ lies over a branch point of $\gamma$,
we see that $P_{2,W}^0$ is generated by the elements of $\rho^\ast \rm{Pic}^0(C)$ 
and $P_{2,W}^{T' - T''}$.
Note that every element of $P_{2,W}^{T' - T''}$ satisfies equation \eqref{involution invariance}.
In addition, the elements of $\rho^\ast \rm{Pic}^0(C)$ satisfying \eqref{involution invariance} 
are the elements of the Prym variety $Prym(C/B)$. Thus, $P_{2,W}^0/\bar{\iota}$
is isomorphic to $m_{i_1}m_{i_2} \dots m_{i_s}$ copies of $Prym(C/B)$.
\end{proof}

\section{The local geometry of a jump}
\label{jump}
We finish this paper with some remarks concerning rank-2 vector bundles on $X$ with jumps.
We begin by fixing some notation. Let $E$ be such a bundle and assume, for simplicity, 
that it has a single jump over the smooth fibre $T = \pi^{-1}(x_0)$. 
The multiplicity $\mu$ of the jump is then less than or equal to the  non-negative integer 
$n_E := -ch_2(E)$ \cite{Brinzanescu-Moraru1} and the spectral cover of the bundle has two components:
\[ S_E = \mu(\{ x_0 \} \times T^\ast) + \overline{C}, \]
where $\overline{C}$ is a bisection of $J(X)$ whose degree over $T^\ast$ is equal to $n_E - \mu$.
Finally, if $\delta$ is the determinant of $E$, the restriction of $E$ to the fibre $T$ is then of the 
form $\lambda \oplus (\lambda^\ast \otimes \delta_{x_0})$, for some 
$\lambda \in \text{Pic}^{-h}(T)$, $h>0$.
The integer $h$ is called the {\it height} of the jump at $T$;
the height and multiplicity satisfy the inequality $h \leq \mu$.

Suppose that $X$ has multiple fibres and consider one of them, say $m_0T_0$. 
As in section \ref{spectral curve}, 
one can associate to $X$ an elliptic quasi-bundle
$\pi': X' \rightarrow B'$, over an $m_0$-cyclic covering $\varepsilon: B' \rightarrow B$,
such that $T'_0 := \psi^{-1}(T_0) \subset X'$
is a smooth fibre of $\pi'$, where $\psi: X' \rightarrow X$ is the
$m_0$-cyclic covering induced by $\varepsilon$.  
We then say that $E$ has a jump over $T_0$
if and only if the restriction of $\psi^\ast E$ to the fibre $T'_0$ is unstable. 
Consequently, the height and multiplicity of the jump of $E$ over $T_0$
are defined as the height and multiplicity of the jump of $\psi^\ast E$ over $T'_0$.

\subsection{The jumping sequence of a vector bundle}
\label{jumping sequence}
A key ingredient in the description of a jump is the notion of elementary
modification of a vector bundle; 
for its basic properties, we refer the reader to \cite{F}.
Consider a rank-2 vector bundle $E$ on $X$ 
with $\det(E) = \delta$, $c_2(E) = c_2$, and spectral cover 
$S_E = (\sum_{i=1}^k \{ x_i \} \times T^\ast ) + \overline{C}$. 
If $N$ is a line bundle on the smooth fibre $T_0 := \pi^{-1}(x_0)$ 
for which there exists at least one surjection $E|_{T_0} \rightarrow N$,  
then any elementary modification $V$ of $E$ determined by $N$ has invariants
$\det(V) = \delta(-T_0)$, $c_2(V) = c_2+\deg(N)$, and spectral cover
\[ S_V = \deg(N)\left( \{ x_0 \}\times T^\ast \right) + 
\left(\sum_{i=1}^k \{ x_i \} \times T^\ast \right) + \overline{C}. \]
Note that the restriction of $V$ to $T_0$ is of the form
$N \oplus (N^\ast \otimes \delta_{x_0})$.

Suppose that $E$ is unstable on the fibre $T$,
splitting as $E|_T = \lambda \oplus (\lambda^\ast \otimes \delta_{x_0})$
for some $\lambda \in \text{Pic}^{-h}(T)$, $h>0$.
Then, up to a multiple of the identity, 
there is a {\em unique} surjection $E|_T \rightarrow \lambda$,
which defines a canonical elementary modification of $E$ that we denote $\bar{E}$. 
This elementary modification is called {\em allowable} \cite{F}.
Consequently, we can associate to $E$ a finite sequence
$\{ \bar{E}_1, \bar{E}_2, \dots , \bar{E}_l \}$
of allowable elementary modifications such that $\bar{E}_l$ is
the only element of the sequence that does not have a jump at $T$.
This sequence induces the following important invariants of a jump.

\begin{definition}\label{jumping sequence - definition}
Let $T$ be a smooth fibre of $\pi$. Suppose that the vector bundle $E$ 
has a jump over $T$ and consider the corresponding sequence of allowable elementary modifications
$\{ \bar{E}_1, \bar{E}_2, \dots , \bar{E}_l \}$. 
The integer $l$ is called
the {\it length} of the jump at $T$. The {\it jumping sequence} of $T$ is defined
as the set of integers $\{ h_0, h_1, \dots , h_{l-1} \}$, where $h_0 = h$
is the height of $E$ and $h_i$ is the height of
$\bar{E}_i$, for $0 < i \leq l-1$.

If the vector bundle $E$ has a jump over a multiple fibre $m_0T_0$ of $\pi$, we define
the length and jumping sequence of $T_0$ to be the length and jumping sequence of 
the jump of $\psi^\ast E$ over the smooth fibre $T'_0 = \psi^{-1}(T_0)$ of $\psi$, 
where $\psi: X' \rightarrow X$ is the $m_0$-cyclic covering defined at the beginning this section.
\end{definition}

\begin{proposition}
Let $E$ be a rank-2 vector bundle on a non-K\"{a}hler elliptic surface that has a jump of
multiplicity $\mu$ at $T$ with jumping sequence $\{ h_0, h_1, \dots , h_{l-1} \}$.
Then,

(i) For all $0 \leq k \leq l-1$, we have $h_{k+1} \leq h_k$.

(ii) Its allowable elementary modification
$\bar{E}$ has a jump of length $l-1$ over $T$ with
jumping sequence $\{ h_1, \dots , h_{l-1} \}$.

(iii) The multiplicity of the jump is the sum of the jumping sequence,
that is, $\mu = \sum_{i=0}^{l-1} h_i$.

(iv)
An elementary modification $\widetilde{E}$ of $E$ 
determined by a line bundle of degree $r \geq h_0$ has a
jump of length $l+1$ over $T$ with jumping sequence $\{ r, h_0, h_1, \dots , h_{l-1} \}$.
Moreover, the allowable elementary modification of $\widetilde{E}$ is $E$.
\hfill \qedsymbol
\end{proposition}

\subsection{Generic jumps} 

If one starts with a vector bundle $E$ that is semistable on every fibre of $\pi$, 
then one can introduce jumps over the smooth fibres of $\pi$
by performing elementary modifications; however,
this implies the existence of a surjection 
from $E$ to a positive line bundle $N$ over such a fibre $T$. 
If the bundle $E$ is regular, then such a surjection always exists.
But if the bundle $E$ is not regular on the fibre $T$, 
these surjections fail to exist precisely for line bundles $N$ of degree one; 
in this case, one cannot add 
a jump of multiplicity one to $E$ over the fibre $T$.

Let us fix a spectral cover $S = (\sum_{i=1}^k \{ x_i \} \times T^\ast ) + \overline{C}$
such that $\pi^{-1}(x_i)$ is a smooth fibre, for all $i = 1, \dots, k$.
Suppose that the bisection $\overline{C}$ is invariant with respect to an involution
of $J(X)$ determined by a line bundle $\delta$ on $X$.
The generic bundle with spectral cover $S$ can be constructed as follows.
Start with a bundle $E_0$ that has determinant $\delta(kT)$ and spectral cover $\overline{C}$,
which is regular on the fibres $\pi^{-1}(x_i)$. If the surface $X$ does not have multiple
fibres, then this is always possible (see Proposition \ref{existence of regular bundles}).
Fix a line bundle $N$ on $T$ of degree one. 
Finally, perform $k$ elementary modifications with respect to $N$ on the fibres $\pi^{-1}(x_i)$, 
counting multiplicity, to obtain a vector bundle whose jumps have the correct multiplicities.
In this case, every jump has jumping sequence $\{ 1, \dots, 1 \}$.
This construction is, however, far from canonical because it depends on the following:
a choice of rank-2 vector bundle $E_0$, a choice of line bundle $N$ of degree one
and surjections to $N$. 

In general, we can assign to a jump $\pi^{-1}(x_i)$ any jumping sequence 
$\{ h_0, \dots, h_{l-1} \}$, as long as
$\sum_{i=1}^{l-1} h_i$ is equal to the multiplicity of the vertical component
$\{ x_i \} \times T^\ast$ in $S$, in which case we will have to choose a different line bundle
$N$ for each distinct integer in $\{ h_0, \dots, h_{l-1} \}$.
For a detailed classification of such jumps, we refer the reader to 
\cite{Moraru, Brinzanescu-Moraru2}.

\end{document}